# Some operator identities related to q-Hermite polynomials

*Johann Cigler*

Fakultät für Mathematik, Universität Wien

johann.cigler@univie.ac.at

**Abstract**

Some $q-$ analogues of the normal ordering of the operator $(X+sD)^n$ on the polynomials are derived.

## 1. Introduction

Consider the multiplication operator $X$ on the polynomials in $x$ defined by $Xf(x) = xf(x)$ and the differentiation operator $D$ which satisfies $Df(x) = f'(x)$.

Define polynomials $H_n(x,s)$ by

$$H_n(x,s) = xH_{n-1}(x,s) + (n-1)sH_{n-2}(x,s) \tag{1.1}$$

with initial values $H_0(x,s)=1, H_1(x,s)=x$. These are a variant of the Hermite polynomials.

Then the following operator identity holds.

**Theorem1 (J.L. Burchnall [1])**

$$(X+sD)^n = \sum_{k=0}^{n} \binom{n}{k} H_{n-k}(X,s) s^k D^k. \tag{1.2}$$

As far as I know formula (1.2) is (in an equivalent form) due to J. L. Burchnall [1]. See also Gian-Carlo Rota [7], p.45. Their proof uses the fact that

$$(X+sD) = e^{-\frac{X^2}{2s}}(sD)e^{\frac{X^2}{2s}}, \tag{1.3}$$

which follows from the operator form

$$Df(X) = f(X)D + f'(X) \tag{1.4}$$

of the differentiation rule for products $D\big(f(x)g(x)\big) = f(x)D\big(g(x)\big) + f'(x)g(x)$.

I shall give another proof which will be generalized in the sequel.

Let $G_n(x,s) = (X+sD)^n 1$.

This means that

$$G_n(x,s) = xG_{n-1}(x,s) + sG'_{n-1}(x,s). \tag{1.5}$$

Then $G_n(x,s)$ is of the form $G_n(x,s) = \sum_{j=0}^{n} c(n,j) s^j x^{n-2j}$. This gives

$c(0,j) = [j=0]$, $c(n,0) = 1$ and

$c(n,j) = c(n-1,j) + (n+1-2j)c(n-1,j-1)$.

From

$$\frac{n!}{2^j j!(n-2j)!} - \frac{(n-1)!}{2^j j!(n-1-2j)!} - \frac{(n-1)!(n+1-2j)}{2^{j-1}(j-1)!(n+1-2j)!} = \frac{(n-1)!}{2^j j!(n-2j)!}(n-(n-2j)-2j) = 0$$

we conclude that $c(n,j) = \dfrac{n!}{2^j j!(n-2j)!}$.

Thus

$$G_n(x,s) = \sum_j s^j \frac{n!}{2^j j!(n-2j)!} x^{n-2j} \tag{1.6}$$

and

$$G'_n(x,s) = nG_{n-1}(x,s). \tag{1.7}$$

Thus

$G_n(x,s) = xG_{n-1}(x,s) + (n-1)sG_{n-2}(x,s)$. Together with the initial values $G_0(x,s) = 1$ and $G_1(x,s) = x$ this implies $G_n(x,s) = H_n(x,s)$.

Therefore we get the following well-known properties

$$(X+sD)^n 1 = H_n(x,s) = \sum_j s^j \frac{n!}{2^j j!(n-2j)!} x^{n-2j} = e^{\frac{sD^2}{2}} x^n \tag{1.8}$$

and

$$H'_n(x,s) = nH_{n-1}(x,s). \tag{1.9}$$

The multiplication operators $H_n(X,s)$ defined by $H_n(X,s)f(x) = H_n(x,s)f(x)$ satisfy

$$DH_n(X,s) = H_n(X,s)D + nH_{n-1}(X,s). \tag{1.10}$$

Now it is easy to derive (1.2). For

$$(X+sD)^{n+1} = (X+sD)\sum_k \binom{n}{k} H_{n-k}(X,s)(sD)^k = \sum_k \binom{n}{k} XH_{n-k}(X,s)(sD)^k$$

$$+ \sum_k \binom{n}{k} H_{n-k}(X,s)(sD)^{k+1} + \sum_k \binom{n}{k}(n-k)H_{n-k-1}(X,s)(sD)^k$$

$$= \sum_k \binom{n}{k}\left(XH_{n-k}(X,s) + (n-k)H_{n-k-1}(X,s)\right)(sD)^k + \sum_k \binom{n}{k-1} H_{n-k+1}(X,s)(sD)^k$$

$$= \sum_k \binom{n}{k} H_{n-k+1}(X,s)(sD)^k + \sum_k \binom{n}{k-1} H_{n-k+1}(X,s)(sD)^k = \sum_{k=0}^n \binom{n+1}{k} H_{n+1-k}(X,s)(sD)^k.$$

Changing $k \to n-k$ in (1.2) and setting $k = m+j$ we get

**Corollary 1**

$$(X+sD)^n = \sum_{m=0}^n \sum_{j=0}^{Min(m,n-m)} \begin{Bmatrix} n \\ m \end{Bmatrix}_j s^{n-m} X^{m-j} D^{n-m-j}. \tag{1.11}$$

Here the so called Weyl binomial coefficients $\begin{Bmatrix} n \\ m \end{Bmatrix}_j$ are given by

$$\begin{Bmatrix} n \\ m \end{Bmatrix}_j = \frac{n!}{2^j j!(m-j)!(n-m-j)!}. \tag{1.12}$$

In this form (1.2) has been rediscovered several times (cf. e.g. A. Varvak [8]) as normal ordering of $(X+sD)^n$, i.e. as a representation in terms of operators of the form $X^k D^\ell$.

Note that

$$\begin{Bmatrix} n \\ m \end{Bmatrix}_j = \begin{Bmatrix} n \\ n-m \end{Bmatrix}_j = \binom{n-2j}{m-j}\begin{Bmatrix} n \\ j \end{Bmatrix}_j \tag{1.13}$$

and

$$H_n(x,s) = \sum_j \frac{n!}{2^j j!(n-2j)!} s^j x^{n-2j} = \sum_{2j \le n} \begin{Bmatrix} n \\ j \end{Bmatrix}_j s^j x^{n-2j}. \tag{1.14}$$

**2. Normal ordering of** $(X + q^{n-1}sD_q)(X + q^{n-2}sD_q)\cdots(X + sD_q).$

I have found two $q$-analogues connected with a variant of the discrete $q$-Hermite polynomials which give simple formulas.

Let

$$h_n(x,s) = \sum_{j=0}^{\lfloor \frac{n}{2} \rfloor} q^{j^2} s^j \frac{[n]!}{(1+q)(1+q^2)\cdots(1+q^j)[j]![n-2j]!} x^{n-2j} \qquad (2.1)$$

be this variant and let $D_q$ be the $q$-differentiation operator defined by

$$D_q f(x) = \frac{f(qx) - f(x)}{(q-1)x}. \qquad (2.2)$$

Then it is well known (cf. e.g. [2]) that

$$h_n(x,s) = (X + q^{n-1}sD_q)(X + q^{n-2}sD_q)\cdots(X + sD_q)1. \qquad (2.3)$$

For a direct proof of this assertion let

$$h_n(x,s) = (X + q^{n-1}sD)(X + q^{n-2}sD)\cdots(X + sD)1 = \sum_{j=0}^{n} c(n,j)s^j x^{n-2j}. \qquad (2.4)$$

This implies $c(0,j) = [j=0]$ and

$c(n,j) = c(n-1, j) + [n+1-2j]q^{n-1}c(n-1, j-1).$

It is now easily verified that $c(n,0) = 1$ and

$$c(n,j) = q^{j^2}\begin{bmatrix} n \\ 2j \end{bmatrix}[2j-1]!! = q^{j^2}\frac{[n]!}{(1+q)(1+q^2)\cdots(1+q^j)[j]![n-2j]!}.$$

For this is true for $n = 0$. With induction we get that these values satisfy

$c(n, j) - c(n-1, j) - [n+1-2j]q^{n-1}c(n-1, j-1)$

$= q^{j^2}\frac{[n]!}{(1+q)(1+q^2)\cdots(1+q^j)[j]![n-2j]!} - q^{j^2}\frac{[n-1]!}{(1+q)(1+q^2)\cdots(1+q^j)[j]![n-1-2j]!}$

$- q^{n-1}[n+1-2j]q^{(j-1)^2}\frac{[n-1]!}{(1+q)(1+q^2)\cdots(1+q^{j-1})[j-1]![n+1-2j]!}$

$= q^{j^2}\frac{[n-1]!}{(1+q)(1+q^2)\cdots(1+q^j)[j]![n-2j]!}\Big([n] - [n-2j] - q^{n-2j}(1+q^j)[j]\Big) = 0.$

Using the $q$-analogue

$$E_q(z) = \sum_{k \geq 0} \frac{q^{\binom{k}{2}} z^k}{[k]_q!} \qquad (2.5)$$

of the exponential series this can be expressed in the form

$$h_n(x,s) = E_{q^2}\left(\frac{qsD_q^2}{[2]_q}\right) x^n. \qquad (2.6)$$

For

$$E_{q^2}\left(\frac{qsD_q^2}{[2]_q}\right) x^n = \sum_j \frac{q^{j^2} s^j}{(1+q)^j [j]_{q^2}!} D_q^{2j} x^n = \sum_j \frac{q^{j^2} s^j}{\prod_{i=0}^{j}(1+q^i)[i]_q} D_q^{2j} x^n = \sum_j \frac{q^{j^2} s^j [n]_q!}{\prod_{i=0}^{j}(1+q^i)[j]_q! [n-2j]_q!} x^{n-2j}.$$

This implies

$$D_q h_n(x,s) = [n] h_{n-1}(x,s). \qquad (2.7)$$

By (2.4) we get the well-known recurrence relation (cf. e.g. [2])

$$h_n(x,s) = x h_{n-1}(x,s) + q^{n-1} s [n-1] h_{n-2}(x,s). \qquad (2.8)$$

For the operator

$$F(n) = \left(X + q^{n-1} sD_q\right)\left(X + q^{n-2} sD_q\right)\cdots\left(X + sD_q\right) \qquad (2.9)$$

we get

**Theorem 2**

$$\left(X + q^{n-1} sD_q\right)\left(X + q^{n-2} sD_q\right)\cdots\left(X + sD_q\right) = \sum_{k=0}^{n} g_n(k, X, s) s^k D_q^k \qquad (2.10)$$

*with*

$$g_n(k,x,s) = \begin{bmatrix} n \\ k \end{bmatrix} \sum_{j=0}^{\lfloor \frac{n-k}{2} \rfloor} s^j q^{j^2 + kj + \binom{k}{2}} \begin{bmatrix} n-k \\ 2j \end{bmatrix} [2j-1]!! \prod_{i=0}^{k-1} \frac{1 + q^{n-j-i}}{1 + q^{j+1+i}} x^{n-k-2j}. \qquad (2.11)$$

*For $k = 0$ this reduces to $g_n(0, x, s) = h_n(x, s)$.*

**Proof**

$$F(n) = \sum_{k=0}^{n} g_n(k,X,s)s^k D_q^k = \left(X + q^{n-1}sD\right)\sum_{k=0}^{n-1} g_{n-1}(k,X,s)s^k D_q^k$$

implies

$$g_n(k,x,s) = x g_{n-1}(k,x,s) + q^{n-1}g_{n-1}(k-1,qx,s) + q^{n-1}s\left(D_q g_{n-1}(k,x,s)\right).$$

Let now

$$g_n(k,x,s) = \sum_j c(n,k,j) s^j x^{n-k-2j}. \tag{2.12}$$

Then

$$c(n,k,j) = c(n-1,k,j) + q^{2n-1-k-2j}c(n-1,k-1,j) + q^{n-1}[n+1-k-2j]c(n-1,k,j-1) \tag{2.13}$$

It now suffices to show that

$$c(n,k,j) = \begin{bmatrix} n \\ k \end{bmatrix} q^{j^2+kj+\binom{k}{2}} \begin{bmatrix} n-k \\ 2j \end{bmatrix}[2j-1]!!\frac{(1+q^{n-j})(1+q^{n-j-1})\cdots(1+q^{n-j-k+1})}{(1+q^{j+1})(1+q^{j+2})\cdots(1+q^{j+k})}$$

$$= \frac{[n]! q^{j^2+kj+\binom{k}{2}}(1+q^{n-j})(1+q^{n-j-1})\cdots(1+q^{n-j-k+1})}{[k]![n-k-2j]![j]!(1+q)(1+q^2)\cdots(1+q^{j+k})}.$$

This is verified by the following computation:

$$c(n,k,j) - c(n-1,k,j) - q^{2n-1-k-2j}c(n-1,k-1,j) - q^{n-1}[n+1-k-2j]c(n-1,k,j-1)$$

$$= \frac{[n-1]! q^{j^2+kj}(1+q^{n-j-1})\cdots(1+q^{n-j-k+1})}{[k]![n-k-2j]![j]!(1+q)(1+q^2)\cdots(1+q^{j+k})}$$

$$\left([n]q^{\binom{k}{2}}(1+q^{n-j}) - q^{\binom{k}{2}}(1+q^{n-j-k})[n-k-2j] - q^{2n-3-3j+\binom{k-2}{2}}[k](1+q^{j+k}) - q^{n-2j+\binom{k-1}{2}}(1+q^{n-j})[j](1+q^{j+k})\right) = 0.$$

**Corollary 2**

$$\left(X + q^{n-1}sD_q\right)\left(X + q^{n-2}sD_q\right)\cdots\left(X + sD_q\right)$$

$$= \sum_{m=0}^{n} \sum_{j=0}^{\text{Min}(m,n-m)} \frac{q^{\binom{j+1}{2}+\binom{n-m}{2}}\left(1+q^{m+1}\right)\left(1+q^{m+2}\right)\cdots\left(1+q^{n-j}\right)[n]! s^{n-m}}{(1+q)(1+q^2)\cdots(1+q^{n-m})[j]![m-j]![n-m-j]!} X^{m-j} D_q^{n-m-j}. \tag{2.14}$$

## 3. Normal ordering of $(X+qsD_q)(X+q^3sD_q)\cdots(X+q^{2n-1}sD_q)$.

Another more natural $q$-analogue of Theorem 1 is

**Theorem 3**

$$(X+qsD_q)(X+q^3sD_q)\cdots(X+q^{2n-1}sD_q) = \sum_{k=0}^{n}\begin{bmatrix}n\\k\end{bmatrix}q^{kn}h_{n-k}(X,s)(sD_q)^k. \qquad (3.1)$$

**Proof**

Consider

$$G(n) = (X+qsD_q)(X+q^3sD_q)\cdots(X+q^{2n-1}sD_q).$$

We show first that $G(n)1 = h_n(x,s)$.

Let

$$f_n(x,s) = (X+sD_q)(X+q^3sD_q)\cdots(X+q^{2n-1}sD_q)1 = \sum_{j=0}^{n}c(n,j)s^j x^{n-2j}. \qquad (3.2)$$

Then $f_n(x,s) = xf_{n-1}(x,q^2 s) + qsD_q f_n(x,q^2 s)$.

This implies $c(0,j) = [j=0]$ and

$$c(n,j) = q^{2j}c(n-1,j) + [n+1-2j]q^{2j-1}c(n-1,j-1).$$

This gives $c(n,0) = 1$ and $c(n,j) = q^{j^2}\begin{bmatrix}n\\2j\end{bmatrix}[2j-1]!! = q^{j^2}\dfrac{[n]!}{(1+q)(1+q^2)\cdots(1+q^j)[j]![n-2j]!}.$

This is true for $n = 0$. With induction we get that these values satisfy

$$c(n,j) - q^{2j}c(n-1,j) - [n+1-2j]q^{2j-1}c(n-1,j-1)$$

$$= q^{j^2}\frac{[n]!}{(1+q)(1+q^2)\cdots(1+q^j)[j]![n-2j]!} - q^{j^2+2j}\frac{[n-1]!}{(1+q)(1+q^2)\cdots(1+q^j)[j]![n-1-2j]!}$$

$$-q^{2j-1}[n+1-2j]q^{(j-1)^2}\frac{[n-1]!}{(1+q)(1+q^2)\cdots(1+q^{j-1})[j-1]![n+1-2j]!}$$

$$= q^{j^2}\frac{[n-1]!}{(1+q)(1+q^2)\cdots(1+q^j)[j]![n-2j]!}\Big([n] - q^{2j}[n-2j] - (1+q^j)[j]\Big) = 0.$$

Therefore we get $f_n(x,s) = h_n(x,s)$.

This implies that $h_n(x,s)$ also satisfies the recurrence

$$h_n(x,s) = xh_{n-1}(x,q^2s) + qsD_q h_n(x,q^2s). \qquad (3.3)$$

From (2.1) it is clear that $h_n(qx,q^2s) = q^n h_n(x,s)$.

The proof of (3.1) follows from

$$(X+qsD_q)\sum_k \begin{bmatrix} 2 \\ k \end{bmatrix} q^{kn}(q^2s)^k h_{2-k}(X,q^2s)D_q^k = \sum_k \begin{bmatrix} n \\ k \end{bmatrix} q^{kn}(q^2s)^k Xh_{n-k}(X,q^2s)D_q^k$$

$$+qs\sum_k \begin{bmatrix} n \\ k \end{bmatrix} q^{kn}(q^2s)^k h_{n-k}(qX,q^2s)D_q^{k+1} + qs\sum_k \begin{bmatrix} n \\ k \end{bmatrix} q^{kn}(q^2s)^k \left(D_q h_{n-k}(X,q^2s)\right)D_q^k$$

$$=\sum_k \begin{bmatrix} n \\ k \end{bmatrix} q^{kn}(q^2s)^k \left(Xh_{n-k}(X,q^2s) + qs\left(D_q h_{n-k}(X,q^2s)\right)\right)D_q^k + \sum_k \begin{bmatrix} n \\ k \end{bmatrix} q^{(k+1)(n+1)} s^k h_{n-k}(X,s)D_q^{k+1}$$

$$=\sum_k q^k \begin{bmatrix} n \\ k \end{bmatrix} q^{k(n+1)} s^k h_{n+1-k}(X,s)D_q^k + \sum_k \begin{bmatrix} n \\ k-1 \end{bmatrix} q^{k(n+1)} s^k h_{n-k+1}(X,s)D_q^k$$

$$=\sum_{k=0}^n \begin{bmatrix} n+1 \\ k \end{bmatrix} q^{k(n+1)} s^k h_{n+1-k}(X,s)D_q^k.$$

**Corollary 3**

$$(X+qsD_q)(X+q^3sD_q)\cdots(X+q^{2n-1}sD_q)$$

$$=\sum_{m=0}^n \sum_{j=0}^{Min(m,n-m)} \frac{q^{n^2+j^2-(m+j)n}[n]!s^{n-m}}{(1+q)(1+q^2)\cdots(1+q^j)[j]![m-j]![n-m-j]!} X^{m-j} D_q^{n-m-j}.$$

**4. Normal ordering of $(X+sD_q)^n$.**

The first terms of $(X+sD_q)^n$ in the ordering of Corollary 1 are

$sD_q + X,$

$s^2 D_q^2 + (1+q)sXD + s + X^2,$

$s^3 D_q^3 + (1+q+q^2)s^2 XD_q^2 + (2+q)s^2 D_q + (1+q+q^2)sX^2 D_q + (2+q)sX + X^3,$

$s^4 D_q^4 + (1+q+q^2+q^3)s^3 XD_q^3 + (3+2q+q^2)s^3 D_q^2 + (1+q+2q^2+q^3+q^4)s^2 X^2 D_q^2$

$+(3+5q+3q^2+q^3)s^2 XD_q + (2+q)s^2 + (1+q+q^2+q^3)sX^3 D_q + (3+2q+q^2)sX^2 + X^4.$

Consider the variant of the $q$-Hermite polynomials introduced in [7] by

$$H_n(x,s\,|\,q) = \left(X + sD_q\right)^n 1. \tag{4.1}$$

They are related to the $q$-Lucas polynomials, which have been studied in [3] and [4]. We define these $q$-Lucas polynomials by

$$L_n(x,s) = \sum_{k=0}^{\lfloor \frac{n}{2} \rfloor} q^{\binom{k}{2}} \frac{[n]}{[n-k]} \begin{bmatrix} n-k \\ k \end{bmatrix} s^k x^{n-2k} \tag{4.2}$$

for $n > 0$ with initial value

$$L_0(x,s) = 1. \tag{4.3}$$

Observe that the choice of initial value is different from the one used in [3].

The first values are $1, x, x^2 + (1+q)s, x^3 + (1+q+q^2)sx, x^4 + (1+q+q^2+q^3)sx^2 + (q+q^3)s^2, \cdots$.

It is easily verified that they satisfy (cf. [3])

$$(X + (1-q)sD_q)L_n(x,-s) = L_{n+1}(x,-s) + sL_{n-1}(x,-s) \tag{4.4}$$

for $n \geq 2$,

$$\left(X + (1-q)sD_q\right)L_1(x,-s) = x^2 + (1-q)s = L_2(x,-s) + sL_0(x,-s) + s \tag{4.5}$$

and

$$\left(X + (1-q)sD_q\right)L_0(x,-s) = x = L_1(x,-s). \tag{4.6}$$

As shown in [5]

$$H_n(x,(1-q)s\,|\,q) = \left(X + (1-q)sD_q\right)^n 1 = \sum_{j=0}^{\lfloor \frac{n}{2} \rfloor} \binom{n}{j} s^j L_{n-2j}(x,-s). \tag{4.7}$$

In order to make this paper self-contained we give a new proof:

(4.7) is obviously true for $n = 0$ and $n = 1$. In the general case we get

$$H_{2n+1}(x,(1-q)s\,|\,q) = \left(X+(1-q)sD_q\right)H_{2n}(x,s\,|\,q) = \left(X+(1-q)sD_q\right)\sum_{j=0}^{n}\binom{2n}{j}s^j L_{2n-2j}(x,-s)$$

$$= \sum_{j=0}^{n}\binom{2n}{j}s^j\left(X+(1-q)sD_q\right)L_{2n-2j}(x,-s) = \sum_{j=0}^{n-1}\binom{2n}{j}s^j\left(L_{2n+1-2j}(x,-s)+sL_{2n-1-2j}(x,-s)\right)+\binom{2n}{n}s^n x$$

$$= \sum_{j=0}^{n-1}\binom{2n}{j}s^j L_{2n+1-2j}(x,-s) + \sum_{j=1}^{n-1}\binom{2n}{j-1}s^j L_{2n+1-2j}(x,-s) + \binom{2n}{n-1}s^n x + \binom{2n}{n}s^n x = \sum_{j=0}^{n}\binom{2n+1}{j}s^j L_{2n+1-2j}(x,-s)$$

and

$$H_{2n}(x,(1-q)s\,|\,q) = \left(X+(1-q)sD_q\right)H_{2n-1}(x,s\,|\,q) = \left(X+(1-q)sD_q\right)\sum_{j=0}^{n-1}\binom{2n-1}{j}s^j L_{2n-1-2j}(x,-s)$$

$$= \sum_{j=0}^{n-1}\binom{2n-1}{j}s^j\left(X+(1-q)sD_q\right)L_{2n-1-2j}(x,-s) = \sum_{j=0}^{n-2}\binom{2n-1}{j}s^j\left(L_{2n-2j}(x,-s)+sL_{2n-2-2j}(x,-s)\right)$$

$$+\binom{2n-1}{n-1}s^{n-1}\left(x^2+(1-q)s\right)$$

$$= \sum_{j=0}^{n-2}\binom{2n-1}{j}s^j L_{2n-2j}(x,-s) + \sum_{j=1}^{n-1}\binom{2n-1}{j-1}s^j L_{2n-2j}(x,-s) + \binom{2n-1}{n-1}s^{n-1}\left(x^2+(1-q)s\right)$$

$$= \sum_{j=0}^{n-2}\binom{2n}{j}s^j L_{2n-2j}(x,-s) + \binom{2n-1}{n-2}s^{n-1}L_2(x,-s) + \binom{2n-1}{n}s^{n-1}(x^2+(1-q)s)$$

$$= \sum_{j=0}^{n-2}\binom{2n}{j}s^j L_{2n-2j}(x,-s) + \binom{2n}{n-1}s^{n-1}L_2(x,-s) + \binom{2n}{n}s^n = \sum_{j=0}^{n}\binom{2n}{j}s^j L_{2n-2j}(x,-s).$$

For the next theorem we need a generalization $L_n^{(k)}(x,s)$ of the Lucas polynomials. We define them by

$$L_n^{(k)}(x,s) = \sum_{j=0}^{\left\lfloor \frac{n}{2}\right\rfloor} q^{\binom{j}{2}} \frac{[n+k]}{[n+k-j]}\begin{bmatrix}n+k-j\\k\end{bmatrix}\begin{bmatrix}n-j\\j\end{bmatrix} s^j x^{n-2j} \tag{4.8}$$

with $L_0^{(k)}(x,s) = 1$. It is clear that $L_n^{(0)}(x,s) = L_n(x,s)$.

Computer calculations led to

**Theorem 4**

$$\left(X + (1-q)sD_q\right)^n = \sum_{k=0}^{n} A(n,k,X)(1-q)^k s^k D_q^k \tag{4.9}$$

with

$$A(n,k,x) = \sum_{i=0}^{\lfloor \frac{n-k}{2} \rfloor} \binom{n}{i} s^i L_{n-2i-k}^{(k)}(x,-s). \tag{4.10}$$

Before proving this theorem let us make some remarks.

$\sum_{k=0}^{n} A(n,k,X)D_q^k$ is a linear combination of terms of the form $x^{n-2i-k-2j}D_q^k$.

Let $m = i+k+j, \ell = i+j$ be the uniquely determined integers such that $2i+k+2j = m+\ell$ and $k = m-\ell$. Then we get

$$A(n,k,x) = \sum_{i=0}^{\lfloor \frac{n-k}{2} \rfloor} \binom{n}{i} s^i \sum_{j=0}^{\lfloor \frac{n-k-2i}{2} \rfloor} q^{\binom{j}{2}} \frac{[n-2i]}{[n-2i-j]} \begin{bmatrix} n-2i-j \\ k \end{bmatrix} \begin{bmatrix} n-k-2i-j \\ j \end{bmatrix} (-s)^j x^{n-k-2i-2j}$$
$$= \sum_{i,m,\ell} \binom{n}{i} (-1)^{\ell-i} s^\ell q^{\binom{\ell-i}{2}} \frac{[n-2i]}{[n-i-\ell]} \begin{bmatrix} n-i-\ell \\ m-\ell \end{bmatrix} \begin{bmatrix} n-m-i \\ \ell-i \end{bmatrix} x^{n-m-\ell} \tag{4.11}$$

Comparing coefficients we see that (4.9) is equivalent with

$$\left(X + sD_q\right)^n = \sum_{m,\ell} \begin{Bmatrix} n \\ m \end{Bmatrix}_{\ell,q} X^{n-m-\ell} s^m D_q^{m-\ell}, \tag{4.12}$$

where the $q$ – Weyl binomial coefficients are given by

$$\begin{Bmatrix} n \\ m \end{Bmatrix}_{\ell,q} = \frac{1}{(1-q)^\ell} \sum_i \binom{n}{i} (-1)^{\ell-i} q^{\binom{\ell-i}{2}} \frac{[n-2i]}{[n-i-\ell]} \begin{bmatrix} n-i-\ell \\ m-\ell \end{bmatrix} \begin{bmatrix} n-m-i \\ \ell-i \end{bmatrix}. \tag{4.13}$$

For the special case $m = \ell$ we get

$$H_n(x, s \mid q) = \sum_\ell \begin{Bmatrix} n \\ \ell \end{Bmatrix}_{\ell, q} x^{n-2\ell} s^\ell \tag{4.14}$$

with

$$\begin{Bmatrix} n \\ \ell \end{Bmatrix}_{\ell, q} = \frac{1}{(1-q)^\ell} \sum_i \binom{n}{i} (-1)^{\ell-i} q^{\binom{\ell-i}{2}} \frac{[n-2i]}{[n-i-\ell]} \begin{bmatrix} n-\ell-i \\ \ell-i \end{bmatrix}. \tag{4.15}$$

Comparing (4.13) with (4.15) we see that

$$\begin{Bmatrix} n \\ m \end{Bmatrix}_{\ell, q} = \begin{bmatrix} n - 2\ell \\ m - \ell \end{bmatrix} \begin{Bmatrix} n \\ \ell \end{Bmatrix}_{\ell, q}. \tag{4.16}$$

My original proof of (4.9) has been rather clumsy. But as has been observed by
J. Zeng [9] Theorem 4 follows immediately from (4.15), which has been proved in [5]
and [6], and (4.16) which has been proved by A. Varvak [8], Theorem 6.4.

The same idea can be used to give a direct proof of Theorem 4:

Define coefficients $\begin{Bmatrix} n \\ m \end{Bmatrix}_{\ell, q}$ by

$$(X + sD_q)^n = \sum_{m, \ell} \begin{Bmatrix} n \\ m \end{Bmatrix}_{\ell, q} X^{m-\ell} s^{n-m} D_q^{n-m-\ell}.$$

Then

$$\sum_{m, \ell} \begin{Bmatrix} n+1 \\ m \end{Bmatrix}_{\ell, q} X^{m-\ell} s^{n+1-m} D_q^{n+1-m-\ell} = (X + sD_q)^{n+1} = (X + sD_q) \sum_{m, \ell} \begin{Bmatrix} n \\ m \end{Bmatrix}_{\ell, q} X^{m-\ell} s^{n-m} D_q^{n-m-\ell}.$$

Since $D_q f(X) = D_q(f(X)) + f(qX) D_q$ we get

$$\begin{Bmatrix} n+1 \\ m \end{Bmatrix}_{\ell, q} = \begin{Bmatrix} n \\ m-1 \end{Bmatrix}_\ell + [m+1-\ell] \begin{Bmatrix} n \\ m \end{Bmatrix}_{\ell-1, q} + q^{m-\ell} \begin{Bmatrix} n \\ m \end{Bmatrix}_{\ell, q}. \tag{4.17}$$

This recurrence together with the initial values $\begin{Bmatrix} 0 \\ 0 \end{Bmatrix}_{0,q} = 1$ and all other $\begin{Bmatrix} 0 \\ m \end{Bmatrix}_{\ell,q} = 0$ determines $\begin{Bmatrix} n \\ m \end{Bmatrix}_{\ell,q}$ uniquely if we set $\begin{Bmatrix} n \\ m \end{Bmatrix}_{\ell,q} = 0$ for $m < 0$ or $\ell < 0$.

It now suffices to show that

$$\begin{Bmatrix} n \\ m \end{Bmatrix}_{\ell,q} := \begin{bmatrix} n - 2\ell \\ m - \ell \end{bmatrix} \begin{Bmatrix} n \\ \ell \end{Bmatrix}_{\ell,q} \tag{4.18}$$

satisfies (4.17) and the initial values.

First observe that (4.15) is an easy consequence of (4.7):

By (4.14) and (4.7) we get

$$H_n(x, (1-q)s \mid q) = \sum_{\ell} \begin{Bmatrix} n \\ \ell \end{Bmatrix}_{\ell,q} x^{n-2\ell}(1-q)^{\ell} s^{\ell} = \sum_{j=0}^{\lfloor \frac{n}{2} \rfloor} \binom{n}{j} s^{j} L_{n-2j}(x, -s)$$

$$= \sum_{i=0}^{\lfloor \frac{n}{2} \rfloor} \binom{n}{i} s^{i} \sum_{k=0}^{\lfloor \frac{n-2i}{2} \rfloor} q^{\binom{k}{2}} \frac{[n-2i]}{[n-2i-k]} \begin{bmatrix} n - 2i - k \\ k \end{bmatrix} (-s)^{k} x^{n-2i-2k}.$$

Setting $i + k = \ell$ and comparing coefficients of $x^{n-2\ell}$ we get

$$\begin{Bmatrix} n \\ \ell \end{Bmatrix}_{\ell,q} (1-q)^{\ell} s^{\ell} = \sum \binom{n}{i} (-1)^{\ell-i} q^{\binom{\ell-i}{2}} \frac{[n-2i]}{[n-i-\ell]} \begin{bmatrix} n - i - \ell \\ \ell - i \end{bmatrix} s^{\ell},$$

which gives (4.15).

Now we must verify (4.17). This becomes

$$\begin{bmatrix} n+1-2\ell \\ m-\ell \end{bmatrix} \begin{Bmatrix} n+1 \\ \ell \end{Bmatrix}_{\ell,q} = \begin{bmatrix} n-2\ell \\ m-1-\ell \end{bmatrix} \begin{Bmatrix} n \\ \ell \end{Bmatrix}_{\ell,q} + [m+1-\ell] \begin{bmatrix} n-2\ell+2 \\ m-\ell+1 \end{bmatrix} \begin{Bmatrix} n \\ \ell-1 \end{Bmatrix}_{\ell-1,q} + q^{m-\ell} \begin{bmatrix} n-2\ell \\ m-\ell \end{bmatrix} \begin{Bmatrix} n \\ \ell \end{Bmatrix}_{\ell,q}$$

$$= \left( \begin{bmatrix} n-2\ell \\ m-1-\ell \end{bmatrix} + q^{m-\ell} \begin{bmatrix} n-2\ell \\ m-\ell \end{bmatrix} \right) \begin{Bmatrix} n \\ \ell \end{Bmatrix}_{\ell,q} + [m+1-\ell] \begin{bmatrix} n-2\ell+2 \\ m-\ell+1 \end{bmatrix} \begin{Bmatrix} n \\ \ell-1 \end{Bmatrix}_{\ell-1,q}$$

$$= \begin{bmatrix} n+1-2\ell \\ m-\ell \end{bmatrix} \begin{Bmatrix} n \\ \ell \end{Bmatrix}_{\ell,q} + [m+1-\ell] \begin{bmatrix} n-2\ell+2 \\ m-\ell+1 \end{bmatrix} \begin{Bmatrix} n \\ \ell-1 \end{Bmatrix}_{\ell-1,q}$$

and reduces to

$$\begin{Bmatrix} n+1 \\ \ell \end{Bmatrix}_{\ell,q} = \begin{Bmatrix} n \\ \ell \end{Bmatrix}_{\ell,q} + [n+2-2\ell] \begin{Bmatrix} n \\ \ell-1 \end{Bmatrix}_{\ell-1,q}$$

which is clear by (4.14) and (4.1).

Thus we get again Theorem 4.

# References

[1] J. L. Burchnall, A note on the polynomials of Hermite, Quart. J. Math. 12(1941), 9-11

[2] J. Cigler, Elementare q-Identitäten, Sém. Lotharingien Comb., B05a (1981), 29 pp.

[3] J. Cigler, A new class of $q-$Fibonacci polynomials, Electr. J. Comb. 10 (2003),#R 19

[4] J. Cigler, $q-$Lucas polynomials and associated Rogers-Ramanujan type identities, arXiv:0907.0165

[5] J. Cigler and J. Zeng, A curious $q-$analogue of Hermite polynomials, arXiv:0905.0228

[6] M. Josuat-Vergès, Rook placements in Young diagrams and permutation enumeration, arXiv: 0811.0524

[7] G.-C. Rota, Finite Operator Calculus, Academic Press 1975

[8] A. Varvak, Rook numbers and the normal ordering problem, J. Comb. Th. A 112 (2005), 292-307

[9] J. Zeng, A remark on the $q-$Weyl binomial coefficients, Personal communication, October 11, 2010